\documentclass{amsart}
\usepackage{amssymb}
\usepackage{latexsym}
\usepackage{amsthm}
\usepackage{amscd}
\theoremstyle{plain}

\newtheorem{thm}{Theorem}[section]

\newtheorem{lem}[thm]{Lemma}
\newtheorem{prop}[thm]{Proposition}

\newtheorem{defn}{Definition}[section]
\newtheorem{rem}{Remark}[section]
\newtheorem{ex}{Example}[section]
\newtheorem{pr}{Problem}

\numberwithin{equation}{section}

\newcommand{\Z}{\mathbb Z}
\newcommand{\C}{\mathbb C}
\newcommand{\Q}{\mathbb Q}
\newcommand{\R}{\mathbb R}
\newcommand{\PP}{\mathbb P}
\newcommand{\N}{\mathbb N}
\newcommand{\HH}{\mathbb H}

\allowdisplaybreaks[1]

\begin{document}

\title[DIFFERENTIAL EQUATIONS SATISFIED BY MODULAR FORMS AND K3 SURFACES]
{Differential equations satisfied by modular forms and K3 surfaces}
\author{Yifan Yang}
\address{Department of Applied Mathematics, National Chiao Tung
University, Hsinchu 300, TAIWAN}
\email{yfyang@math.nctu.edu.tw}
\author{Noriko Yui}
\address{Department of Mathematics and Statistics, Queen's University,
Kingston, Ontario Canada K7L 3N6}
\email{yui@mast.queensu.ca}
\date{\today}
\thanks{N. Yui was partially supported by Discovery Grant from the
  Natural Science and Engineering Research Council (NSERC) of 
  Canada. Y. Yang was supported by Grant 93-2115-M-009-014 from the
  National Science Council (NSC) of Taiwan and by the National Center
  for Theoretical Sciences (NCTS) of Taiwan.}


\begin{abstract}
  We study differential equations satisfied by modular forms 
  of two variables associated to $\Gamma_1\times \Gamma_2$, where
  $\Gamma_i$ ($i=1,2$) are genus zero subgroups of $SL_2(\R)$ commensurable
  with $SL_2(\Z)$, e.g., $\Gamma_0(N)$ or $\Gamma_0(N)^*$ for some $N$. 
  In some examples, these differential equations are realized as the
  Picard--Fuchs differential equations of families of $K3$ surfaces with
  large Picard numbers, e.g., $19, 18, 17, 16$.  Our method rediscovers
  some of the Lian--Yau examples of ``modular relations'' involving power
  series solutions to the second and the third order differential equations
  of Fuchsian type in \cite{LY1, LY2}. 
\end{abstract}

\maketitle

\section{Introduction}
\label{sect1}

Lian and Yau \cite{LY1,LY2} studied arithmetic properties of mirror
maps of pencils of certain $K3$ surfaces, and further, they considered
mirror maps of certain families of Calabi--Yau threefolds \cite{LY3}.
Lian and Yau observed in a number of explicit examples a mysterious
relationship (now the so-called {\it mirror moonshine phenomenon})
between mirror maps and the McKay--Thompson series (Hauptmoduls of one
variable associated to a genus zero congruence subgroup of $SL_2(\R)$)
arising from the Monster. Inspired by the work of Lian and
Yau, Verrill--Yui \cite{VY} further computed more examples of mirror
maps of one-parameter families of lattice polarized $K3$ surfaces with
Picard number $19$. The outcome of Verrill--Yui's calculations
suggested that the mirror maps themselves are not always Hauptmoduls,
but they are commensurable with Hauptmoduls (referred as the
modularity of mirror maps). This fact was indeed established 
by Doran \cite{Dor} for $M_n$-lattice polarized $K3$ surfaces of 
Picard number $19$ with maximal unipotent monodromy
(where $M_n=U\perp (-E_8)^2\perp\langle-2n\rangle$).
More generally, Doran \cite{Dor2} considered the commensurability of
``maximal $n$-dimensional families of rank $20-n$ lattice polarized
families of $K3$ surfaces, and he showed that all such families of
$K3$ surfaces are commensurable to autormorphic forms. 
 
The mirror maps were calculated via the Picard--Fuchs differential equations of the $K3$
families in question. Therefore, the determination of the
Picard--Fuchs differential equations played the central role in
their investigations.
  
In this paper, we will address the inverse problem of a kind.  That is,
instead of starting with families of $K3$ surfaces or families of
Calabi--Yau threefolds, we start with modular forms and functions of
more than one variable associated to certain subgroups of $SL_2(\R)$. 

More specifically, the main focus our discussions in this paper are on
modular forms and functions of two variables.  Here is the precise
definition.
\smallskip

\begin{defn} {\rm
  Let $\HH$ denote the upper half-plane $\{\tau:\Im\tau>0\}$, and let
  $\HH^\ast=\HH\cup\Q\cup\{\infty\}$. Let $\Gamma_1$ and $\Gamma_2$ be
  two subgroups of $SL_2(\R)$ commensurable with $SL_2(\Z)$. We call a
  function $F:\HH^\ast\times\HH^\ast\longmapsto\C$ of two variables a
  {\it modular form (of two variables) of weight $(k_1,k_2)$ on
  $\Gamma_1\times\Gamma_2$ with character $\chi$} if $F$ is
  meromorphic on $\HH^\ast\times\HH^\ast$ such that
  $$
    F(\gamma_1\tau_1,\gamma_2\tau_2)=\chi(\gamma_1,\gamma_2)
    (c_1\tau_1+d_1)^{k_1}(c_2\tau_2+d_2)^{k_2}F(\tau_1,\tau_2)
  $$
  for all
  $$
    \gamma_1=\begin{pmatrix}a_1&b_1\\ c_1&d_1\end{pmatrix}\in\Gamma_1,
    \qquad
    \gamma_2=\begin{pmatrix}a_2&b_2\\ c_2&d_2\end{pmatrix}\in\Gamma_2.
  $$
  If $F$ is a modular form (of two variables) of weight $(0,0)$ with trivial
  character, then we also call $F$ a {\it modular function (of
  two variables) on $\Gamma_1\times\Gamma_2$.}}
\end{defn}

\smallskip
\noindent{\bf Notation.} We let $q_1=e^{2\pi i\tau_1}$ and
  $q_2=e^{2\pi i\tau_2}$. For a variable $t$ we let $D_t$ denote the
  the differential operator $t\frac{\partial}{\partial t}$.
\smallskip

\begin{rem}\label{remark 1.1}
{\rm Stienstra and Zagier \cite{SZ} have introduced the notion of
{\it bi-modular forms} (of two variables). Let $\Gamma\subset
\mbox{SL}_2(\R)$,  and let $\tau_1,\,\tau_2\in \HH$. Let $k_1,\, k_2$
be integers. A two-variable meromorphic function $F: \HH\times
\HH\to\C$ is called a {\it bi-modular} form of weight $(k_1,k_2)$ on
$\Gamma$ if for any $\gamma=\begin{pmatrix} a & b \\ c & d
\end{pmatrix}\in\Gamma$, it satisfies the transformation formula:
$$
  F(\gamma\tau_1,\gamma\tau_2)=(c\tau_1+d)^{k_1}(c\tau_2+d)^{k_2}
  F(\tau_1,\tau_2).
$$
For instance,
$$
  F(\tau_1, \tau_2)=\tau_1-\tau_2
$$
is a bi-modular form for $\mbox{SL}_2(\Z)$ of weight $(-1,-1)$.
Another typical example is
$$
  F(\tau_1,\tau_2)=E_2(\tau_1)-\frac{1}{\tau_1-\tau_2},
$$
which is a bi-modular form of weight $(2,0)$ for $\mbox{SL}_2(\Z)$.

For bi-modular forms of Stienstra--Zagier, the fundamental domain
$\HH\times\HH/\Gamma$ is not of finite volume.
On the other hand, for our modular forms (of two variables),
the fundamental domain is $\HH/\Gamma_1\times \HH/\Gamma_2$, which
is always of finite volume.

We should emphasize that the two notions of two variable modular forms
(namely, our modular forms and bi-modular forms of Stienstra and Zagier) 
are indeed different. Also we mention that our modular
forms are not a special case of Hilbert modular forms.}
\end{rem}
\smallskip

  The problems that we will consider here are formulated
as follows : {\it Given a modular form $F$ (of two variables), 
determine a differential equation it satisfies, and
construct a family of $K3$ surfaces (or degenerations of a family of
Calabi--Yau threefolds at some limit points) having the determined
differential equation as its Picard--Fuchs differential equation.} 
This kind of problem may be called a {\it geometric realization 
problem}.
 
In fact, a similar problem was already considered by Lian 
and Yau in their papers \cite{LY1, LY2}. They discussed the so-called
``modular relations'' involving power series solutions to second and
third order differential equations of Fuchsian type (e.g., hypergeometric
differential equations $_2F_1,\, _3F_2$) and modular forms of weight 
$4$ using mirror symmetry. More recently, van Enckevort and van Straten
\cite{ES05} considered the following geometric realization problem: 
{\it Starting with a certain forth order
differential equation whose monodromy representation can be calculated,
find a one-parameter families of Calabi--Yau threefolds (if it exists), whose associated
Picard--Fuchs differential equation is the given one}. Also a recent article
of Doran and Morgan \cite{DorMor05} addressed the geometric realization
question in the context of an old question of Griffiths: {\it When does an
integral variation of Hodge structure come from geometry?}. A rigorous
answer was presented for one-parameter
families of Calabi--Yau threefolds with $h^{2,1}=1$ with generalized
Picard--Fuch differential eqations, relating mirror symmetry
and integral variations of Hodge structure. 
\smallskip
 
In this paper, we will focus our discussion on modular forms (of
two variables) of weight $(1,1)$.  We will determine the differential 
equations satisfied by modular forms (of two variables) of weight $(1,1)$ 
associated to $\Gamma_1\times \Gamma_2$ where $\Gamma_i$ are genus zero 
subgroups of $SL_2(\R)$ of the form $\Gamma_0(N)$ and $\Gamma_0(N)^*$. Then 
the existence and the construction of particular modular forms of weight 
$(1,1)$ are discussed, using solutions of some hypergeometric differential 
equations.  Moreover, we determine the differential equations they satisfy.  
  Further, several examples of modular forms (of two variables) and their
  differential equations are discussed aiming to realize these differential
  equations as the Picard--Fuchs differential equations
  of some families of $K3$ surfaces (or degenerations of families of 
  Calabi--Yau threefolds) with large Picard numbers $19,18,17$
  and $16$.

  It should be pointed out that our paper and our results have non-empty
  intersections with the results of Lian and Yau \cite{LY1, LY2}. Indeed,
  our approach rediscovers some of the examples of Lian and Yau.

   Our contributions may be summarized as follows. From geometric
   point of veiw, we give examples of two-parameter families of $K3$
   surfaces which after pull-back along a morphism from
   $(t_1,t_2)$-space to $(x,y)$-space decouple as a direct product of
   two one-parameter families of elliptic curves. From function
   theoretic point of veiw, we give examples of non-trivial
   substitions transforming certain (two-variables) GKZ hypergeometic
   hypergeometric functions into a product of two (one-variable) GKZ
   hypergeometric fucntions.
Finally, from moduli point of view, we give examples of moduli spaces for
$K3$ surfaces with extra structure and show that these moduli spaces
are quotients of ${\mathfrak H}\times{\mathfrak H}$.
 
\section{Differential equations satisfied by modular forms (of
two variables)}
\label{sect2}

  We will now determine differential equations satisfied by
  modular forms (of two variables) of weight $(1,1)$ on $\Gamma_1\times
\Gamma_2$.
\smallskip

\begin{thm} \label{DE satisfied by modular forms}
 {\sl  Let $F(\tau_1,\tau_2)$ be a modular form (of two variables) of 
weight $(1,1)$, and let $x(\tau_1,\tau_2)$ and $y(\tau_1,\tau_2)$ be 
non-constant modular functions (of two variables) on 
$\Gamma_1\times\Gamma_2$, where $\Gamma_i$ ($i=1,2$) are
subgroups of $SL_2(\R)$ commensurable with $SL_2(\Z)$. Then $F$, as a
function of $x$ and $y$, satisfy a system of partial differential
equations 
  \begin{equation} \label{main}
  \begin{split}
    D_x^2F+a_0D_xD_yF+a_1D_xF+a_2D_yF+a_3F=0, \\
    D_y^2F+b_0D_xD_yF+b_1D_xF+b_2D_yF+b_3F=0,
  \end{split}
  \end{equation}
  where $a_i$ and $b_i$ are algebraic functions of $x$ and $y$, and
  can be expressed explicitly as follows. Suppose that, for each
  function $t$ among $F$, $x$, and $y$, we let
  $$
    G_{t,1}=\frac{D_{q_1}t}t=\frac1{2\pi i}\frac{dt}{t\,d\tau_1}, \qquad
    G_{t,2}=\frac{D_{q_2}t}t=\frac1{2\pi i}\frac{dt}{t\,d\tau_2}.
  $$
  Then we have
  $$
    a_0=\frac{2G_{y,1}G_{y,2}}{G_{x,1}G_{y,2}+G_{y,1}G_{x,2}}, \qquad
    b_0=\frac{2G_{x,1}G_{x,2}}{G_{x,1}G_{y,2}+G_{y,1}G_{x,2}},
  $$
  $$
    a_1=\frac{G_{y,2}^2(D_{q_1}G_{x,1}-2G_{F,1}G_{x,1})
      -G_{y,1}^2(D_{q_2}G_{x,2}-2G_{F,2}G_{x,2})}
      {G_{x,1}^2G_{y,2}^2-G_{y,1}^2G_{x,2}^2},
  $$
  $$
    b_1=\frac{-G_{x,2}^2(D_{q_1}G_{x,1}-2G_{F,1}G_{x,1})
      +G_{x,1}^2(D_{q_2}G_{x,2}-2G_{F,2}G_{x,2})}
      {G_{x,1}^2G_{y,2}^2-G_{y,1}^2G_{x,2}^2},
  $$
  $$
    a_2=\frac{G_{y,2}^2(D_{q_1}G_{y,1}-2G_{F,1}G_{y,1})
      -G_{y,1}^2(D_{q_2}G_{y,2}-2G_{F,2}G_{y,2})}
      {G_{x,1}^2G_{y,2}^2-G_{y,1}^2G_{x,2}^2},
  $$
  $$
    b_2=\frac{-G_{x,2}^2(D_{q_1}G_{y,1}-2G_{F,1}G_{y,1})
      +G_{x,1}^2(D_{q_2}G_{y,2}-2G_{F,2}G_{y,2})}
      {G_{x,1}^2G_{y,2}^2-G_{y,1}^2G_{x,2}^2},
  $$
  $$
    a_3=-\frac{G_{y,2}^2(D_{q_1}G_{F,1}-G_{F,1}^2)
      -G_{y,1}^2(D_{q_2}G_{F,2}-G_{F,2}^2)}
      {G_{x,1}^2G_{y,2}^2-G_{y,1}^2G_{x,2}^2},
  $$
  and
  $$
    b_3=-\frac{-G_{x,2}^2(D_{q_1}G_{F,1}-G_{F,1}^2)
      +G_{x,1}^2(D_{q_2}G_{F,2}-G_{F,2}^2)}
      {G_{x,1}^2G_{y,2}^2-G_{y,1}^2G_{x,2}^2}.
  $$}
\end{thm}

In order to prove Theorem \ref{DE satisfied by modular forms}, we first
need the following lemma, which is an analogue of the classical
Ramanujan's differential equations
$$D_qE_2=\frac{E_2^2-E_4}{12}=-24\sum_{n\in\N}\frac{n^2q^n}{(1-q^n)^2},$$
$$ D_qE_4=\frac{E_2E_4-E_6}3=240\sum_{n\in\N}\frac{n^4q^n}{(1-q^n)^2},$$
$$ D_qE_6=\frac{E_2E_6-E_4^2}2=\sum_{n\in\N}\frac{n^6q^n}{(1-q^n)^2}$$
where 
\begin{equation} \label{Ek}
  E_k=1-\frac{2k}{B_k}\sum_{n\in\N}\frac{n^{k-1}q^n}{1-q^n}
\end{equation}
are the Eisenstein series of weight $k$ on $SL_2(\Z)$,
where $B_k$ denotes the $k$-th Bernoulli number, e.g.,
$B_2=\frac{1}{6},\, B_4=-\frac{1}{30}$ and $B_6=\frac{1}{42}$.
\smallskip

\begin{lem} \label{lemma 2.2} {\sl We retain the notations of Theorem
  \ref{DE satisfied by modular forms}. Then 

  (a) $G_{x,1}$ and $G_{y,1}$ are modular forms (of two variables) 
of weight $(2,0)$, 

  (b) $G_{x,2}$ and $G_{y,2}$ are modular forms (of two variables)
of weight $(0,2)$, 

  (c) $D_{q_1}G_{x,1}-2G_{F,1}G_{x,1}$, $D_{q_1}G_{y,1}-2G_{F,1}G_{y,1}$ and 
   $D_{q_1}G_{F,1}-G_{F,1}^2$ are modular forms (of two variables) 
of weight $(4,0)$, and 

  (d) $D_{q_2}G_{x,2}-2G_{F,2}G_{x,2}$, $D_{q_2}G_{y,2}-2G_{F,2}G_{y,2}$
and $D_{q_2}G_{F,2}-G_{F,2}^2$ are modular forms (of two
variables) of weight $(0,4)$.}   
\end{lem}
\begin{proof} We shall prove (a) and (c); the proof of (b) and (d) is
  similar.

  By assumption, $x$ is a modular function (of two variables) on
  $\Gamma_1\times\Gamma_2$. That is, for all
  $\gamma_1=\begin{pmatrix}a_1&b_1\\ c_1&d_1\end{pmatrix}\in\Gamma_1$
  and all $\gamma_2=\begin{pmatrix}a_2&b_2\\
  c_2&d_2\end{pmatrix}\in\Gamma_2$, one has
  $$
    x(\gamma_1\tau_1,\gamma_2\tau_2)=x(\tau_1,\tau_2)
  $$
  Taking the logarithmic derivatives of the above equation with
  respect to $\tau_1$, we obtain
  $$
    \frac1{(c_1\tau_1+d_1)^2}\frac{\dot x}x(\gamma_1\tau_1,\tau_2)
   =\frac{\dot x}x(\tau_1,\tau_2),
  $$
  or
  \begin{equation} \label{temp1: lemma 2.2}
    G_{x,1}(\gamma_1\tau_1,\gamma_2\tau_2)=(c_1\tau_1+d_1)^2
    G_{x,1}(\tau_1,\tau_2),
  \end{equation}
  where we let $\dot x$ denote the derivative of the two-variable
  function $x$ with respect to the first variable. This shows that
  $G_{x,1}$ is a modular form of weight $(2,0)$ on
  $\Gamma_1\times\Gamma_2$ with the trivial character. The proof for
  the case $G_{y,1}$ is similar.

  Likewise, taking the logarithmetic derivatives of the equation
  $$
    F(\gamma_1\tau_1,\gamma_2\tau_2)=\chi(\gamma_1,\gamma_2)
    (c_1\tau_1+d_1)(c_2\tau_2+d_2)F(\tau_1,\tau_2)
  $$
  with respect to $\tau_1$, we obtain
  $$
    \frac1{(c_1\tau_1+d_1)^2}\frac{\dot F}F(\gamma_1\tau_1,\gamma_2\tau_2)
   =\frac{c_1}{(c_1\tau_1+d_1)}+\frac{\dot F}F(\tau_1,\tau_2),
  $$
  or, equivalently
  \begin{equation} \label{temp2: lemma 2.2}
    G_{F,1}(\gamma_1\tau_1,\gamma_2\tau_2)=\frac{c_1(c_1\tau_1+d_1)}
    {2\pi i}+(c_1\tau_1+d_1)^2G_{F,1}(\tau_1,\tau_2).
  \end{equation}
  Now, differentiating (\ref{temp1: lemma 2.2}) with respect to
  $\tau_1$ again, we obtain
  $$
    \frac{\dot G_{x,1}}{(c_1\tau_1+d_1)^2}
   (\gamma_1\tau_1,\gamma_2\tau_2)
   =2c_1(c_1\tau_1+d_1)G_{x,1}(\tau_1,\tau_2)
   +(c_1\tau_1+d_1)^2\dot G_{x,1}(\tau_1,\tau_2),
  $$
  or
  $$
    D_{q_1}G_{x,1}(\gamma_1\tau_1,\gamma_2\tau_2)
   =\frac{c_1(c_1\tau_1+d_1)^3}{\pi i}G_{x,1}(\tau_1,\tau_2)
   +(c_1\tau_1+d_1)^4D_{q_1}G_{x,1}(\tau_1,\tau_2).
  $$
  On the other hand, we also have, by (\ref{temp1: lemma 2.2}) and
  (\ref{temp2: lemma 2.2}),
  $$
    G_{F,1}G_{x,1}(\gamma_1\tau_1,\gamma_2\tau_2)
    =\frac{c_1(c_1\tau_1+d_1)^3}{2\pi i}G_{x,1}(\tau_1,\tau_2)
    +(c_1\tau_1+d_1)^4G_{F,1}G_{x,1}(\tau_1,\tau_2).
  $$
  From these two equations we see that
  $D_{q_1}G_{x,1}-2G_{F,1}G_{x,1}$ is a modular form (of two 
variables) of weight $(4,0)$ with the trivial character.

  Finally, differentiating (\ref{temp2: lemma 2.2}) with respect to
  $\tau_1$ and multiplying by $(c_1\tau_1+d_1)^2$ we have
  \begin{equation*}
  \begin{split}
    D_{q_1}G_{F,1}(\gamma_1\tau_1,\gamma_2\tau_2)
   &=\frac{c_1^2(c_1\tau_1+d_1)^2}{(2\pi i)^2}
   +\frac{c_1(c_1\tau_1+d_1)^3}{\pi i}G_{F,1}(\tau_1,\tau_2) \\
   &\qquad\qquad+(c_1\tau_1+d_1)^4D_{q_1}G_{F,1}(\tau_1,\tau_2).
  \end{split}
  \end{equation*}
  Combining this with the square of (\ref{temp2: lemma 2.2}) we see
  that $D_{q_1}G_{F,1}-G_{F,1}^2$ is a modular form  
  of weight $(4,0)$ on $\Gamma_1\times\Gamma_2$. This completes the 
  proof of the lemma.
\end{proof}
\smallskip

\begin{proof}[Proof of Theorem \ref{DE satisfied by modular forms}] In
  light of Lemma \ref{lemma 2.2}, the functions $a_k$, $b_k$ are all
  modular functions on $\Gamma_1\times\Gamma_2$, and thus can be
  expressed as algebraic functions of $x$ and $y$. Therefore, it
  suffices to verify (\ref{main}) as formal identities.  By the chain
  rule we have
  $$
    \begin{pmatrix} D_{q_1}F \\ D_{q_2}F\end{pmatrix}
   =\begin{pmatrix} x^{-1}D_{q_1}x & y^{-1}D_{q_1}y \\
      x^{-1}D_{q_2}x & y^{-1}D_{q_2}y \end{pmatrix}
    \begin{pmatrix} D_x F \\ D_y F \end{pmatrix}.
  $$
  It follows that
  $$
    \begin{pmatrix} D_x F \\ D_y F \end{pmatrix}
   =\frac F{G_{x,1}G_{y,2}-G_{x,2}G_{y,1}}
    \begin{pmatrix} G_{y,2} & -G_{y,1} \\
      -G_{x,2} & G_{x,1} \end{pmatrix}
    \begin{pmatrix} G_{F,1} \\ G_{F,2} \end{pmatrix}
  $$
  Writing
  $$
    \Delta=G_{x,1}G_{y,2}-G_{x,2}G_{y,1},
  $$
  and
  $$
    \Delta_x=G_{F,1}G_{y,2}-G_{F,2}G_{y,1}, \qquad
    \Delta_y=-G_{x,2}G_{F,1}+G_{x,1}G_{F,2},
  $$
  we have
  \begin{equation} \label{DxF, DyF}
    D_xF=F\frac{\Delta_x}\Delta, \qquad D_yF=F\frac{\Delta_y}\Delta.
  \end{equation}
  Applying the same procedure on $D_xF$ again, we obtain
  \begin{equation*}
  \begin{split}
    \begin{pmatrix} D_x^2 F \\ D_yD_x F \end{pmatrix}
  &=\frac 1\Delta
    \begin{pmatrix} G_{y,2} & -G_{y,1} \\
      -G_{x,2} & G_{x,1} \end{pmatrix}
    \begin{pmatrix} D_{q_1}(F\Delta_x/\Delta) \\
      D_{q_2}(F\Delta_x/\Delta) \end{pmatrix} \\
  &=\frac F\Delta
    \begin{pmatrix} G_{y,2} & -G_{y,1} \\
      -G_{x,2} & G_{x,1} \end{pmatrix}
    \left\{\frac{\Delta_x}\Delta
    \begin{pmatrix} G_{F,1} \\ G_{F,2} \end{pmatrix}+
    \begin{pmatrix} D_{q_1}(\Delta_x/\Delta) \\
      D_{q_2}(\Delta_x/\Delta) \end{pmatrix}\right\}.
  \end{split}
  \end{equation*}
  That is,
  \begin{equation} \label{DxxF}
    D_x^2F=F\frac{\Delta_x^2}{\Delta^2}+\frac F\Delta\left(
    G_{y,2}D_{q_1}\frac{\Delta_x}\Delta-G_{y,1}D_{q_2}
    \frac{\Delta_x}\Delta\right)
  \end{equation}
  and
  \begin{equation} \label{DyxF}
    D_yD_xF=F\frac{\Delta_x\Delta_y}{\Delta^2}+\frac F\Delta
    \left(-G_{x,2}D_{q_1}\frac{\Delta_x}\Delta+G_{x,1}D_{q_2}
    \frac{\Delta_x}\Delta\right).
  \end{equation}
  We then substitute (\ref{DxF, DyF}), (\ref{DxxF}), and (\ref{DyxF})
  into (\ref{main}) and find that (\ref{main}) indeed holds.
  (The details are tedious, but straightforward calculations. We
  omit the details here.)
\end{proof}

\section{Modular forms (of two variables) associated to solutions of 
hypergeometric differential equations}
\label{sect4}

Here we will construct modular forms (of two variables) of weight $(1,1)$ using
solutions of some hypergeometric differential equations. Our main
result of this section is the following theorem.  

\begin{thm} \label{Hypergeometric DE} {\sl Let $0<a<1$ be a positive
  real number. Let $f(t)=\,_2F_1(a,a;1;t)$ be a solution of the
  hypergeometric differential equation
  \begin{equation} \label{HG DE}
    t(1-t)f^{\prime\prime}+[1-(1+2a)t]f^\prime-a^2f=0.
  \end{equation}
  Let
  $$
    F(t_1,t_2)=f(t_1)f(t_2)(1-t_1)^a(1-t_2)^a,
  $$
  $$
    x=\frac{t_1+t_2}{(t_1-1)(t_2-1)},
    \quad y=\frac{t_1t_2}{(t_1+t_2)^2}.
  $$
  Then $F$ is a modular form of weight $(1,1)$ for $\Gamma_1\times
  \Gamma_2$, provided that $t_1$ and $t_2$ are modular functions (of
one variable) for $\Gamma_1$ and $\Gamma_2$, respectively.
Furthermore, $F$, as a function of $x$ and $y$, 
% but without the condition of modular functions on $t_1$ and $t_2$, 
is a solution of the partial differential equations
  \begin{equation} \label{HGDE x}
    D_x(D_x-2D_y)F+x(D_x+a)(D_x+1-a)F=0,
  \end{equation}
  and
  \begin{equation} \label{HGDE y}
    D_y^2F-y(2D_y-D_x+1)(2D_y-D_x)F=0,
  \end{equation}
  where $D_x=\partial/\partial x$ and $D_y=\partial/\partial y$.  }
\end{thm}
\smallskip

\begin{rem}\label{remark 3.1} {\rm Theorem 2.1 of Lian and Yau \cite{LY2}
is essentially the same as our Theorem \ref{Hypergeometric DE}, though
the formulation and proof are different.

The condition that $t_1,\, t_2$ are modular functions (of one
variable) for $\Gamma_1$ and $\Gamma_2$ are used to draw the conclustion 
that $F$ is a modular form (of two variables) for $\Gamma_1\times\Gamma_2$.
However, the modular property of $t_1,\,t_2$ is irrelevant to derive 
\ref{HGDE x} and \ref{HGDE y} from \ref{HG DE}.}
\end{rem}
\smallskip

We will present our proof of Theorem \ref{Hypergeometric DE} now. 
For this, we need one more ingredient, namely, the Schwarzian derivatives.
\smallskip

\begin{lem} \label{Schwarzian} {\sl Let $f(t)$ and $f_1(t)$ be two
  linearly independent solutions of a differential equation
  $$
    f^{\prime\prime}+p_1f^\prime+p_2f=0.
  $$
  Set $\tau:=f_1(t)/f(t)$. Then the associated
  Schwarzian differential equation
  $$
    2Q\left(\frac{dt}{d\tau}\right)^2+\{t,\tau\}=0,
  $$
  where $\{t,\tau\}$ is the Schwarzian derivative
  $$
    \{t,\tau\}=\frac{dt^3/d\tau^3}{dt/d\tau}-\frac32
    \left(\frac{dt^2/d\tau^2}{dt/d\tau}\right)^2,
  $$
  satisfies
  $$
    Q=\frac{4p_2-2p_1^\prime-p_1^2}4.
  $$}
\end{lem}
\begin{proof}
This is standard, and proof can be found, for instance, in Lian
and Yau \cite{LY3}.
\end{proof}
\smallskip

\begin{proof}[Proof of Theorem \ref{Hypergeometric DE}] Let $f_1$ be
  another solution of (\ref{HG DE}) linearly independent of $f$, and
  set $\tau=f_1/f$. Then a classical identity asserts that
  $$
    f^2=c\exp\left\{-\int^t\frac{1-(1+2a)u}{u(1-u)}\,du\right\}
    \frac{dt}{d\tau}
    =\frac{cdt/d\tau}{t(1-t)^{2a}},
  $$
  where $c$ is a constant depending on the choice of $f_1$. Thus,
  letting 
  $$
    q_1=e^{2\pi if_1(t_1)/f(t_1)}\quad\mbox{and}\quad
    q_2=e^{2\pi if_1(t_2)/f(t_2)},
  $$ 
  the function $F$, with a suitable choice of $f_1$, is in fact
  $$
    F(t_1,t_2)=\left(\frac{D_{q_1}t_1\cdot D_{q_2}t_2}
    {t_1t_2}\right)^{1/2}.
  $$
  We now apply the differential identities in (\ref{main}), which hold
  for arbitrary $F$, $x$, and $y$.
  We have
  $$
    G_{x,1}:=\frac{D_{q_1}x}x=\frac{(1+t_2)D_{q_1}t_1}{(t_1+t_2)(1-t_1)},
    \quad
    G_{x,2}:=\frac{D_{q_2}x}x=\frac{(1+t_1)D_{q_2}t_2}{(t_1+t_2)(1-t_2)},
  $$
  $$
    G_{y,1}:=\frac{D_{q_1}y}y=\frac{(t_2-t_1)D_{q_1}t_1}{t_1(t_1+t_2)},
    \qquad
    G_{y,2}:=\frac{D_{q_2}y}y=\frac{(t_1-t_2)D_{q_2}t_2}{t_2(t_1+t_2)},
  $$
  $$
    G_{F,1}:=\frac{D_{q_1}F}F=\frac{t_1D_{q_1}^2t_1-(D_{q_1}t_1)^2}
      {2t_1D_{q_1}t_1}, \quad
    G_{F,2}:=\frac{D_{q_2}F}F=\frac{t_2D_{q_2}^2t_2-(D_{q_2}t_2)^2}
      {2t_2D_{q_2}t_2}.
  $$
  It follows that
  $$
    a_0:=\frac{2G_{y,1}G_{y,2}}{G_{x,1}G_{y,2}+G_{y,1}G_{x,2}}
    =-\frac{2(t_1-1)(t_2-1)}{t_1t_2+1}=-\frac2{1+x},
  $$
  $$
    b_0:=\frac{2G_{x,1}G_{x,2}}{G_{x,1}G_{y,2}+G_{y,1}G_{x,2}}
    =\frac{2t_1t_2(t_1+1)(t_2+1)}{(t_1-t_2)^2(t_1t_2+1)}
    =\frac{2y(1+2x)}{(1+x)(1-4y)},
  $$
  \begin{equation*}
  \begin{split}
    a_1:&=\frac{G_{y,2}^2(D_{q_1}G_{x,1}-2G_{F,1}G_{x,1})
      -G_{y,1}^2(D_{q_2}G_{x,2}-2G_{F,2}G_{x,2})}
      {G_{x,1}^2G_{y,2}^2-G_{y,1}^2G_{x,2}^2} \\
    &=\frac{t_1+t_2}{t_1t_2+1}=\frac x{1+x},
  \end{split}
  \end{equation*}
  \begin{equation*}
  \begin{split}
    b_1:&=\frac{-G_{x,2}^2(D_{q_1}G_{x,1}-2G_{F,1}G_{x,1})
      +G_{x,1}^2(D_{q_2}G_{x,2}-2G_{F,2}G_{x,2})}
      {G_{x,1}^2G_{y,2}^2-G_{y,1}^2G_{x,2}^2} \\
    &=\frac{t_1t_2(t_1+1)(t_2+1)}{(t_1-t_2)^2(t_1t_2+1)}
     =\frac{y(1+2x)}{(1+x)(1-4y)},
  \end{split}
  \end{equation*}
  $$
    a_2:=\frac{G_{y,2}^2(D_{q_1}G_{y,1}-2G_{F,1}G_{y,1})
      -G_{y,1}^2(D_{q_2}G_{y,2}-2G_{F,2}G_{y,2})}
      {G_{x,1}^2G_{y,2}^2-G_{y,1}^2G_{x,2}^2}=0,
  $$
  \begin{equation*}
  \begin{split}
    b_2:&=\frac{-G_{x,2}^2(D_{q_1}G_{y,1}-2G_{F,1}G_{y,1})
      +G_{x,1}^2(D_{q_2}G_{y,2}-2G_{F,2}G_{y,2})}
      {G_{x,1}^2G_{y,2}^2-G_{y,1}^2G_{x,2}^2} \\
    &=-\frac{2t_1t_2}{(t_1-t_2)^2}=-\frac{2y}{1-4y}.
  \end{split}
  \end{equation*}
Moreover, we have
  \begin{equation*}
  \begin{split}
    a_3:&=-\frac{G_{y,2}^2(D_{q_1}G_{F,1}-G_{F,1}^2)
      -G_{y,1}^2(D_{q_2}G_{F,2}-G_{F,2}^2)}
      {G_{x,1}^2G_{y,2}^2-G_{y,1}^2G_{x,2}^2} \\
    &=\frac{(t_1-1)(t_2-1)(t_1+t_2)\left\{t_1^2\dot t_2^4
     (2\dot t_1\dddot t_1-3\ddot t_1^2)-t_2^2\dot t_1^4
     (2\dot t_2\dddot t_2-3\ddot t_2^2)\right\}}
      {4(t_1-t_2)(t_1^2t_2^2-1)\dot t_1^4\dot t_2^4},
  \end{split}
  \end{equation*}
where, for brevity, we let $\dot t_j$, $\ddot t_j$, $\dddot t_j$
denote the derivatives $D_{q_j}t_j$, $D_{q_j}^2t_j$, and
$D_{q_j}^3t_j$, respectively. To express $a_3$ in terms of $x$ and
$y$, we note that, by Lemma \ref{Schwarzian},
\begin{equation*}
\begin{split}
  2\dot t_j\dddot t_j-3\ddot t_j^2&=-\dot t_j^4\left(
   \frac{-4a^2}{t_j(1-t_j)}-2\frac d{dt_j}\frac{1-(1+2a)t_j}{t_j(1-t_j)}
  -\frac{(1-(1+2a)t_j)^2}{t_j^2(1-t_j)^2}\right) \\
  &=-\frac{(t_j-1)^2+4a(1-a)t_j}{t_j^2(t_j-1)^2}\dot t_j^4.
\end{split}
\end{equation*}
It follows that
$$
  a_3=a(1-a)\frac{t_1+t_2}{t_1t_2+1}=\frac{a(1-a)x}{1+x}.
$$
Likewise, we have
\begin{equation*}
\begin{split}
  b_3:&=-\frac{-G_{x,2}^2(D_{q_1}G_{F,1}-G_{F,1}^2)
    +G_{x,1}^2(D_{q_2}G_{F,2}-G_{F,2}^2)}
    {G_{x,1}^2G_{y,2}^2-G_{y,1}^2G_{x,2}^2} \\
  &=a(1-a)\frac{t_1t_2(t_1+t_2)}{(t_1-t_2)^2(t_1t_2+1)}
   =\frac{a(1-a)xy}{(1+x)(1-4y)}.
\end{split}
\end{equation*}
Then, by (\ref{main}), the function $F$, as a function of $x$ and $y$,
satisfies
\begin{equation} \label{HG 1}
  D_x^2F-\frac2{1+x}D_xD_yF+\frac x{1+x}D_xF+\frac{a(1-a)x}{1+x}F=0
\end{equation}
and
\begin{equation} \label{HG 2}
\begin{split}
 &D_y^2F+\frac{2y(1+2x)}{(1+x)(1-4y)}D_xD_yF
 +\frac{y(1+2x)}{(1+x)(1-4y)}D_xF \\
 &\qquad\qquad -\frac{2y}{1-4y}D_yF+\frac{a(1-a)xy}{(1+x)(1-4y)}F=0.
\end{split}
\end{equation}
Finally, we can deduce the claimed differential equations by taking
(\ref{HG 1}) times $(1+x)$ and (\ref{HG 2}) times $(1-4y)$ minus
(\ref{HG 1}) times $y$, respectively.
\end{proof}

\section{Examples}

\begin{ex}\label{Example 4.1}
{\rm Let $j$ be the elliptic modular $j$-function, and let
$E_4(\tau)=1+240\sum_{n\in\N}\frac{n^3q^n}{1-q^n}$, $q=e^{2\pi i\tau}$, 
be the Eisenstein series of weight $4$ on $SL_2(\Z)$. Set
$$
  x=2\frac{1/j(\tau_1)+1/j(\tau_2)-1728/(j(\tau_1)j(\tau_2))}{1+
    \sqrt{(1-1728/j(\tau_1))(1-1728/j(\tau_2))}}, \qquad
  y=\frac1{j(\tau_1)j(\tau_2)x^2},
$$
and $$F=(E_4(\tau_1)E_4(\tau_2))^{1/4}.$$ Then $F$
% is a modular form
%(of two variables) of weight $(1,1)$ for $SL_2(\Z)\times SL_2(\Z)$,
%and it
satisfies the system of partial differential equations:
$$
  (1-432x)D_x^2F-2D_xD_yF-432xD_x-60xF=0,
$$
$$
  (1-4y)D_y^2F+4yD_xD_yF-yD_x^2F-yD_xF-2yD_yF=0.
$$}

We should remark that the functions $x$ and $y$ are modular functions
(of two variables) for $\Gamma_1\times \Gamma_2$ where 
$\Gamma_1=\Gamma_2$ is a subgroup of $SL_2(\Z)$ of index $2$.  On the
other hand, in the sense of Stienstra-Zagier, $x$ and $y$ are
bi-modular functions for the group $SL_2(\Z)$ (cf. Remark 1.1). 
\end{ex}

We have noticed that this system of differential equation 
belongs to a general class of partial differential equations which involve 
solutions of hypergeometric hypergeometric differential equations
discussed in Theorem \ref{Hypergeometric DE}.

Here we will prove the assertion of Example \ref{Example 4.1} using 
Theorem \ref{Hypergeometric DE}. 

\begin{proof}[Proof of Example 4.1] We first make a change of variable
  $x\mapsto-\bar x/432$. For convenience, we shall denote the new
  variable $\bar x$ by $x$ again. Thus, we are required to show that
  the functions 
  $$
    x=-864\frac{1/j(\tau_1)+1/j(\tau_2)-1728/(j(\tau_1)j(\tau_2))}{1+
      \sqrt{(1-1728/j(\tau_1))(1-1728/j(\tau_2))}}, \qquad
    y=\frac{432^2}{j(\tau_1)j(\tau_2)x^2},
  $$
  and $F=(E_4(\tau_1)E_4(\tau_2))^{1/4}$ satisfy
  $$
  (1+x)D_x^2F-2D_xD_yF+xD_x+\frac 5{36}xF=0,
  $$
  and
  $$
  (1-4y)D_y^2F+4yD_xD_yF-yD_x^2F-yD_xF-2yD_yF=0.
  $$

  For brevity, we let $j_1$ denote $j(\tau_1)$ and $j_2$ denote
  $j(\tau_2)$. We now observe that the function $x$ can be
  alternatively expressed as
\begin{equation*}
\begin{split}
  x&=-864\frac{1/j_1+1/j_2-1728/(j_1j_2)}{1-(1-1728/j_1)(1-1728/j_2)}
     \left(1-\sqrt{(1-1728/j_1)(1-1728/j_2)}\right) \\
   &=\frac12\left(\sqrt{(1-1728/j_1)(1-1728/j_2)}-1\right).
\end{split}
\end{equation*}
Setting
$$
  t_1=\frac{\sqrt{1-1728/j_1}-1}{\sqrt{1-1728/j_1}+1}, \qquad
  t_2=\frac{\sqrt{1-1728/j_2}-1}{\sqrt{1-1728/j_2}+1},
$$
we have
$$
  x=\frac{t_1+t_2}{(t_1-1)(t_2-1)}.
$$
Moreover, the functions $j_k$, written in terms of $t_k$, are
$j_k=432(t_k-1)^2/t_k$ for $k=1,2$. It follows that
$$
  y=\frac{432^2}{j_1j_2x^2}=\frac{t_1t_2}{(t_1+t_2)^2}.
$$
In view of Theorem \ref{Hypergeometric DE}, setting
$$
  t=\frac{\sqrt{1-1728/j(\tau)}-1}{\sqrt{1-1728/j(\tau)}+1}
$$
it remains to show that the function $f(t)=E_4(\tau)^{1/4}(1-t)^{-1/6}$ is a
solution of the hypergeometric differential equation
$$
  t(1-t)f^{\prime\prime}+(1+4t/3)f^\prime-\frac1{36}f=0,
$$
or equivalently, that
$$
  \frac{E_4(\tau)^{1/4}}{(1-t)^{1/6}}=\, _2F_1(1/6,1/6;1;t).
$$
This, however, follows from the classical identity
$$
  E_4(\tau)^{1/4}=\,_2F_1\left(\frac1{12},\frac5{12};1;\frac{1728}{j(\tau)}\right)
$$
and Kummer's transformation formula
\begin{equation*}
\begin{split}
  &\left(\frac{1+\sqrt{1-z}}2\right)^{2a}\,
  _2F_1\left(a,b;a+b+\frac12;z\right) \\
  &\qquad\qquad=\,
  _2F_1\left(2a,a-b+\frac12;a+b+\frac12;\frac{\sqrt{1-z}-1}
  {\sqrt{1-z}+1}\right).
\end{split}
\end{equation*}
This completes the proof of Example \ref{Example 4.1}.
\end{proof}

\begin{rem} 
{\rm The functions $x$ and $y$ in Example \ref{Example 4.1} (up to 
constant multiple) have also appeared in the paper of Lian and 
Yau \cite{LY3}, Corollary 1.2, as the mirror map of the family of 
$K3$ surfaces defined by degree $12$ hypersurfaces in the 
weighted projective space $\PP^3[1,1,4,6]$.  Further, this $K3$
family is derived from the square of a family of elliptic curves in
the weighted projective space $\PP^2[1,2,3]$. (The geometry
behind this phenomenon is the so-called Shoida--Inose structures,
which has been studied in detail by Long \cite{Long} for one-parameter
families of $K3$ surfaces, and their Picard--Fuchs differential 
equations.) Lian and Yau \cite{LY3} proved that the mirror map of 
the $K3$ family can be given in terms of the
elliptic $j$-function, and indeed, by the functions $x$ and $y$ (up to
constant multiple).  We will discuss more examples of families 
of $K3$ surfaces, their Picard--Fuchs differential equations and 
mirror maps in the section 6.}
\end{rem}

Along the same vein, we obtain more examples of modular forms
of weight $(1,1)$ and modular functions on $\Gamma_0(N)\times
\Gamma_0(N)$ for $N=2,3,4$.
\smallskip

\begin{thm}\label{theorem 4.1} 
{\sl We retain the notations of Theorem \ref{HG DE}. Then
the solutions of the differential equations
(\ref{HGDE x}) and (\ref{HGDE y}) for the cases $a=1/2,1/3,1/4,1/6$ can
be expressed in terms of modular forms and modular functions
on $\Gamma_0(N)\times \Gamma_0(N)$ for some $N$. 

(a) For $a=1/2$, they are given by
$$
  F(\tau_1,\tau_2)=\theta_4(\tau_1)^2\theta_4(\tau_2)^2, \qquad
  t=\theta_2(\tau)^4/\theta_3(\tau)^4,
$$
which are modular on $\Gamma_0(4)\times\Gamma_0(4)$.

(b) For $a=1/3$, they are
$$
  F(\tau_1,\tau_2)=\frac12(3E_2(3\tau_1)-E_2(\tau_1))^{1/2}
    (3E_2(3\tau_2)-E_2(\tau_2))^{1/2}, \quad
  t=-27\frac{\eta(3\tau)^{12}}{\eta(\tau)^{12}},
$$
which are modular on $\Gamma_0(3)\times\Gamma_0(3)$.

(c) For $a=1/4$, they are
$$
  F(\tau_1,\tau_2)=(2E_2(2\tau_1)-E_2(\tau_1))^{1/2}(2E_2(2\tau_2)-E_2(\tau_2))^{1/2},
  \quad t=-64\frac{\eta(2\tau)^{24}}{\eta(\tau)^{24}},
$$
which are modular are $\Gamma_0(2)\times\Gamma_0(2)$.

(d) For $a=1/6$, they are given as in Example \ref{Example 4.1}.
\medskip

%Here $$E_2(q)=1-24\sum_{n\in\N}\frac{nq^n}{1-q^n}$$ is the Eisenstein
%series of weight $2$ for $SL_2(\Z)$, 
Here $$\eta(\tau)=q^{1/24}\prod_{n\in\N}(1-q^n), \qquad q=e^{2\pi i\tau}$$
is the Dedekind eta-function, and
$$
\theta_2(\tau)=q^{1/4}\sum_{n\in\Z} q^{n(n+1)},\quad
\theta_3(\tau)=\sum_{n\in\Z} q^{n^2},\quad
\theta_4(\tau)=\sum_{n\in\Z} (-1)^nq^{n^2}$$
are theta-series.}
\end{thm}

\begin{lem} \label{lemma 4.2} Let $\Gamma$ be a subgroup of $SL_2(\R)$
  commensurable with $SL_2(\Z)$. Let $f(\tau)$ be a modular form (of
one variable) of weight $1$, and $t(\tau)$ be a non-constant modular function 
(of one variable) on $\Gamma$. Then, setting
  $$
    G_t=\frac{D_q t}t, \qquad G_f=\frac{D_q f}f,
  $$
  we have
  $$
    D_t^2f+\frac{D_qG_t-2G_fG_t}{G_t^2}D_tf-
    \frac{D_qG_f-G_f^2}{G_t^2}f=0.
  $$
\end{lem}

\begin{proof}[Proof of Theorem 4.1] To prove part (a) we use the
  well-known identities
$$
  \theta_3^2=\,_2F_1\left(\frac12,\frac12;1;
  \frac{\theta_2^4}{\theta_3^4}\right)
$$
(see \cite{Yg04} for a proof using Lemma \ref{lemma 4.2}) and
$$
  \theta_3^4=\theta_2^4+\theta_4^4.
$$
Applying Theorem 3.1 and observing that
$$
  \theta_3^2\left(1-\frac{\theta_2^4}{\theta_3^4}\right)^{1/2}
 =\theta_3^2\frac{\theta_4^2}{\theta_3^2}=\theta_4^2,
$$
we thus obtain the claimed differential equation.

For parts (b), we need to show that the function
$$
  f(\tau)=\frac{(3E_2(3\tau)-E_2(\tau))^{1/2}}{(1-t)^{1/3}}
$$
satisfies
$$
  t(1-t)\frac{d^2}{dt^2}f+(1-5t/3)\frac d{dt}f-\frac19f=0,
$$
or, equivalently,
\begin{equation} \label{temp1: theorem 4.1}
  (1-t)D_t^2f-\frac23tD_tf-\frac19tf=0.
\end{equation}
Let $G_t$ and $G_f$ be defined as in Lemma \ref{lemma 4.2}. For
convenience we also let $g=(3E_2(3\tau)-E_2(\tau))/2$. We have
$$
  G_t=\frac12(3E_2(3\tau)-E_2(\tau))=g
$$
and
$$
  G_f=\frac{D_q g}{2g}-\frac1{3(1-t)}D_q t=
  \frac{D_q g}{2g}+\frac t{3(1-t)}g.
$$
It follows that
$$
  \frac{D_qG_t-2G_fG_t}{G_t^2}=g^{-2}
  \left(D_q g-2\left(\frac{D_q g}{2g}+\frac t{3(1-t)}g\right)g\right)
 =-\frac{2t}{3(1-t)}.
$$
Moreover, we can show that $(D_qG_f-G_f^2)/G_t^2$ is equal to
$-t/(9(1-t))$ by comparing enough Fourier coefficients. This
establishes (\ref{temp1: theorem 4.1}) and hence part (b).

The proof of part (c) is similar, and we shall skip the details here.
\end{proof}

\section{More examples}
\label{sect5}

We may also consider groups like $\Gamma_0(N)^*\times
\Gamma_0(N)^*$ where $\Gamma_0(N)^*$ denotes the group
generated by $\Gamma_0(N)$ and the Atkin--Lehner
involution $w_N=\begin{pmatrix} 0 & -1 \\ N & 0\end{pmatrix}$
for some $N$.  (Note that $\Gamma_0(N)^*$ is contained
in the normalizer of $\Gamma_0(N)$ in $SL_2(\R)$.) 
Also the entire list of $N$ giving rise to genus zero
groups $\Gamma_0(N)^*$ is known (cf.\cite{CMS04}), and we will be
interested in some of those genuz zero groups. We can determine differential
equations satisfied by modular forms (of two variables) of weight $(1,1)$
on $\Gamma_0(N)^*\times \Gamma_0(N)^*$ for some $N$ (giving rise to
genus zero subgroups $\Gamma_0(N)^*$). 

We first prove a generalization of Theorem \ref{Hypergeometric DE}.

\begin{thm}\label{theorem 5.1} Let $0<a,b<1$ be positive real numbers. 
  Let $f(t)=\,_2F_1(a,b;1;t)$ be a solution of the hypergeometric
  differential equation
  \begin{equation} \label{temp: Theorem 5.1}
    t(1-t)f^{\prime\prime}+[1-(1+a+b)t]f^\prime-abf=0.
  \end{equation}
  Set
  $$
    F(t_1,t_2)=f(t_1)f(t_2)(1-t_1)^{(a+b)/2}(1-t_2)^{(a+b)/2},
  $$
  $$
    x=t_1+t_2-2, \qquad y=(1-t_1)(1-t_2).
  $$
  Then $F$, as a function of $x$ and $y$, satisfies

\begin{equation}\label{5.2}
    D_x^2F+2D_xD_yF-\frac1{x+y+1}D_xF+\frac{x}{x+y+1}D_yF
    +\frac{(2ab-a-b)x}{2(x+y+1)}F=0
\end{equation}
  and
  \begin{equation}\label{5.3}
  \begin{split}
   &D_y^2F+\frac{2y}{x^2}D_xD_yF+\frac{y^2}{x^2(x+y+1)}D_xF
    +\frac{y-x-x^2}{x(x+y+1)}D_yF \\
   &\qquad\qquad-\frac{(a+b)(a+b-2)(x^2+x)+(a-b)^2xy-(4ab-2a-2b)y}
  {4x(x+y+1)}F=0.
  \end{split}
  \end{equation}
\end{thm}

\begin{proof} The proof is very similar to that of Theorem
  \ref{Hypergeometric DE}. Let $f_1$ be another solution of the
  hypergeometric differential equation (\ref{temp: Theorem 5.1}), and
  set $\tau:=f_1/f$. We find
  $$
    f^2=c\exp\left\{-\int^t\frac{1-(1+a+b)u}{u(1-u)}\,du\right\}
    \frac{dt}{d\tau}
    =\frac{cdt/d\tau}{t(1-t)^{a+b}}
  $$
  for some constant $c$ depending on the choice of $f_1$. Thus, setting
  $$
    q_1=e^{2\pi if_1(t_1)/f(t_1)}\quad\mbox{and}\quad
    q_2=e^{2\pi if_1(t_2)/f(t_2)},
  $$ 
  we have
  $$
    F(t_1,t_2)=c^\prime\left(\frac{D_{q_1}t_1\cdot D_{q_2}t_2}
    {t_1t_2}\right)^{1/2}
  $$
  for some constant $c^\prime$. We now apply the differential identities
  (\ref{main}). We have, for $j=1,2$,
  $$
    G_{x,j}:=\frac{D_{q_j}x}x=\frac{D_{q_j}t_j}{t_1+t_2-2}, \qquad
    G_{y,j}:=\frac{D_{q_j}y}y=-\frac{D_{q_j}t_j}{1-t_j},
  $$
  and
  $$
    G_{F,j}:=\frac{D_{q_j}F}F=\frac{t_jD_{q_j}^2t_j-(D_{q_j}t_j)^2}
    {2t_jD_{q_j}t_j}.
  $$
  It follows that the coefficients in (\ref{main}) are
  $$
    a_0=2, \qquad
%:=\frac{2G_{y,1}G_{y,2}}{G_{x,1}G_{y,2}+G_{y,1}G_{x,2}}=2,
    b_0
%:=\frac{2G_{x,1}G_{x,2}}{G_{x,1}G_{y,2}+G_{y,1}G_{x,2}}
    =\frac{2(1-t_1)(1-t_2)}{(t_1+t_2-2)^2}=\frac{2y}{x^2},
  $$
  $$
    a_1
    =-\frac1{t_1t_2}=-\frac1{x+y+1}, \qquad
    b_1
     =\frac{(1-t_1)^2(1-t_2)^2}{t_1t_2(t_1+t_2-2)^2}
     =\frac{y^2}{x^2(x+y+1)},
  $$
  $$
    a_2
     =\frac{t_1+t_2-2}{t_1t_2}=\frac x{x+y+1},
  $$
  $$
    b_2
     =-\frac{t_1^2+t_1t_2+t_2^2-2t_1-2t_2+1}{t_1t_2(t_1+t_2-2)}
     =\frac{y-x-x^2}{x(x+y+1)}.
  $$
  Moreover, we have
  \begin{equation*}
  \begin{split}
    a_3
    &=\left\{-\frac{(1-t_1)^2(2\dot t_1\dddot t_1
     -3\ddot t_1^2)}{4(t_1-t_2)\dot t_1^4}
     +\frac{(1-t_2)^2(2\dot t_2\dddot t_2-3\ddot t_2^2)}
      {4(t_1-t_2)\dot t_2^4}
     -\frac{2t_1t_2-t_1-t_2}{4t_1^2t_2^2}\right\} \\
    &\qquad\times(t_1+t_2-2),
  \end{split}
  \end{equation*}
  where we, as before, employ the notations $\dot t_j$, $\ddot t_j$, $\dddot t_j$
  for the derivatives $D_{q_j}t_j$, $D_{q_j}^2t_j$, and
  $D_{q_j}^3t_j$, respectively. Now, by Lemma \ref{Schwarzian}, we have
  $$
    2\dot t_j\dddot t_j-3\ddot t_j^2=\dot t_j^4
    \frac{(a-b)^2t_j^2-(1-t_j)^2+(4ab-2a-2b)t_j}{t_j^2(1-t_j)^2}.
  $$
  It follows that
  $$
    a_3=\frac{(2ab-a-b)(t_1+t_2-2)}{2t_1t_2}=\frac{(2ab-a-b)x}{x+y+1}.
  $$
  A similar calculation shows that
  $$
    b_3=-\frac{(a+b)(a+b-2)(x^2+x)+(a-b)^2xy-(4ab-2a-2b)y}{4x(x+y+1)}.
  $$
  This proves the claimed result.
\end{proof}

\begin{rem}
It should be pointed out that the first identity in our proof of Theorem \ref{theorem 5.1}
is equivalent to the formula in Proposition 4.4 of Lian and Yau \cite{LY1}.
\end{rem}

We now obtain new examples of modular forms of weight $(1,1)$
on $\Gamma_0(N)^*\times \Gamma_0(N)^*$ for some $N$.

\begin{thm} When the pairs of numbers $(a,b)$ in Theorem \ref{theorem
  5.1} are given by $(1/12,5/12)$, $(1/12,7/12)$, $(1/8,3/8)$,
    $(1/8,5/8)$, $(1/6,1/3)$, $(1/6,2/3)$, $(1/4,1/4)$ and
    $(1/4,3/4)$, the solutions $F(t_1,t_2)$ of the differential equations \eqref{5.2}
    and \eqref{5.3} are
    modular forms of weight $(1,1)$ on $\Gamma_0(N)^\ast\times\Gamma_0(N)^\ast$ 
    with $N=1,1,2,2,3,3,4,4$, respectively.
\end{thm}

\begin{proof} We shall prove only the cases $(a,b)=(1/6,1/3)$ and
    $(1/6,2/3)$; the other cases can be proved in the same manner.

  Let
  $$
    s(\tau)=-27\frac{\eta(3\tau)^{12}}{\eta(\tau)^{12}}, \qquad
    E_2(\tau)=1-24\sum_{n=1}^\infty\frac{nq^n}{1-q^n}.
  $$
  From the proof of Part (b) of Theorem \ref{theorem 4.1} we know that
  $$
    f(\tau)=\frac{(3E_2(3\tau)-E_2(\tau))^{1/2}}{(1-s)^{1/3}},
  $$
  as a function of $s$, is equal to $\sqrt2\,_2F_1(1/3,1/3;1;s)$. Now,
  applying the quadratic transformation formula
  $$
    _2F_1(\alpha,\beta;\alpha-\beta+1;x)=(1-x)^{-\alpha}\,_2F_1
    \left(\frac\alpha2,\frac{1+\alpha}2-\beta;\alpha-\beta+1; 
    -\frac{4x}{(1-x)^2}\right)
  $$
  for hypergeometric functions (see, for example \cite[Theorem
  3.1.1]{AAR}) with $\alpha=\beta=1/3$, we obtain
  $$
    (3E_2(3\tau)-E_2(\tau))^{1/2}=\sqrt2\,_2F_1\left(\frac16,\frac13;1;
    -\frac{4s}{(1-s)^2}\right).
  $$
  Observing that the action of the Atkin-Lehner involution $w_3$
  sends $s$ to $1/s$, we find that the function $s/(1-s)^2$ is modular
  on $\Gamma_0(3)^\ast$.  This proves that $F(t_1,t_2)$ is a modular
  form of weight $(1,1)$ for $\Gamma_0(3)^\ast\times\Gamma_0(3)^\ast$ in
  the case $(a,b)=(1/6,1/3)$. 

  Furthermore, an application of another hypergeometric function identity
  $$
    _2F_1(\alpha,\beta;\gamma;x)=(1-x)^{-\alpha}\,_2F_1\left(\alpha,
    \gamma-\beta;\gamma;\frac x{x-1}\right)
  $$
  yields
  $$
    (3E_2(3\tau)-E_2(\tau))^{1/2}=\sqrt2\left(\frac{1-s}{1+s}\right)^{1/3}
    \,_2F_1\left(\frac16,\frac23;1;\frac{4s}{(1+s)^2}\right).
  $$
  This corresponds to the case $(a,b)=(1/6,2/3)$. Again, the function
  $4s/(1+s)^2$ is modular on $\Gamma_0(3)^\ast$.
  This implies that $F(t_1,t_2)$ is a modular form of
  weight $(1,1)$ for $\Gamma_0(3)^\ast\times \Gamma_0(3)^\ast$
  for the case $(a,b)=(1/6,2/3)$.
\end{proof}

\begin{rem} For the remaining pairs $(a,b)$ in Theorem 5.2, we simply list 
the exact expressions of $F(t_1,t_2)$ in terms of modular forms as proofs are similar. 

For $(a,b)=(1/12,5/12)$ and $(1/12,7/12)$, they are
  $$
    \left(\frac{E_6(\tau_1)E_6(\tau_2)}{E_4(\tau_1)E_4(\tau_2)}\right)^{1/2},
    \qquad and \quad
    \left(\frac{E_8(\tau_1)E_8(\tau_2)}{E_6(\tau_1)E_6(\tau_2)}\right)^{1/2},
  $$
  respectively, where $E_k$ are the Eisenstein series in
  (\ref{Ek}). 

For $(a,b)=(1/8,3/8)$ and $(1/8,5/8)$, they are
  $$
    \prod_{j=1}^2\left(\frac{1+s_j}{1-s_j}
   (2E_2(2\tau_j)-E_2(\tau_j))\right)^{1/2}, \quad and \quad
    \prod_{j=1}^2\left(\frac{1-s_j}{1+s_j}
   (2E_2(2\tau_j)-E_2(\tau_j))\right)^{1/2},
  $$
  respectively, where $s_j=-64\eta(\tau_j)^{24}/\eta(\tau_j)^{24}$. 

For $(a,b)=(1/6,1/3)$ and $(1/6,2/3)$, they are
  $$
    \prod_{j=1}^2\left(\frac{1+s_j}{1-s_j}
   (3E_2(3\tau_j)-E_2(\tau_j))\right)^{1/2}, \quad and \quad
    \prod_{j=1}^2\left(\frac{1-s_j}{1+s_j}
   (3E_2(3\tau_j)-E_2(\tau_j))\right)^{1/2},
  $$
  respectively, where $s_j=-27\eta(3\tau_j)^{12}/\eta(\tau_j)$. 

For $(a,b)=(1/4,1/4)$ and $(1/4,3/4)$, they are
  $$
    \prod_{j=1}^2\left(2E_2(2\tau_j)-E_2(\tau_j)\right)^{1/2}, \qquad and \quad
    \prod_{j=1}^2\left(2E_2(2\tau_j)-E_2(\tau_j)\right)^{1/2}
    \frac{1-s_j}{1+s_j},
  $$
  respectively, where $s_j=\theta_2(\tau_j)^4/\theta_3(\tau_j)^4$.
\end{rem}

\bigskip

\section{Picard--Fuchs differential equations of Familes of $K3$ surfaces 
: Part I}
\label{sect6}

One of the motivations of our investigation is to understand the
mirror maps of families of $K3$ surfaces with large Picard nubmers, e.g.,
$19, 18, 17$ or $16$.
Some examples of such families of $K3$ surfaces were discussed
in Lian--Yau \cite{LY2}, Hosono--Lian--Yau \cite{HLY} and 
also in Verrill-Yui \cite{VY}. Some of $K3$ families occured considering
degenerations of Calabi--Yau families.

Our goal here is to construct families of $K3$ surfaces
whose Picard--Fuchs differential equations are given by the
differential equations satisfied by modular forms (of two variables) we
constructed in the earlier sections.  In this section,
we will look into the families of $K3$ surfaces appeared
in Lian and Yau \cite{LY1,LY2}. 
\smallskip

Let $S$ be a $K3$ surface.
We recall some general theory about $K3$ surfaces which are
relevant to our discussion.  We know that
$$H^2(S,\Z)\simeq (-E_8)^2\perp U^3$$
where $U$ is the hyperbolic plane $\begin{pmatrix} 0 & 1 \\ 1 & 0
\end{pmatrix}$
and $E_8$ is the even unimodular negative definite lattice of rank $8$.
The Picard group of $S$, $\mbox{Pic}(S)$, is the group of linear
equivalence classes of Cartier divisors on $S$. Then $\mbox{Pic}(S)$
injects to $H^2(X,\Z)$, and the image of $\mbox{Pic}(S)$ is the algebraic
cycles in $H^2(S,\Z)$.  As $Pic(S)$ is torsion-free, it may be
regarded as a lattice in $H^2(S,\Z)$, called the {\it Picard lattice},
and its rank is denoted by $\rho(S)$.  

According to Arnold--Dolgachev \cite{D96}, two $K3$ surfaces form
a mirror pair $(S, \hat S)$ if
$$\mbox{Pic}(S)^{\perp}_{H^2(S,\Z)}=\mbox{Pic}(\hat S)\perp U\quad
\mbox{as lattices}$$
In terms of ranks, a mirror pair $(S, \hat S)$ is related by
the identity: 
$$22-\rho(S)=\rho(\hat S)+2\Leftrightarrow \rho(S)+\rho(\hat S)=20.$$

\begin{ex}\label{Example 6.1}
{\rm  We will be interested in mirror pairs of
$K3$ surfaces $(S,\hat S)$ whose Picard lattices are of the form
$$Pic(S)=U\quad{and}\quad Pic(\hat S)=U_2\perp (-E_8)^2.$$

We go back to our Example \ref{Example 4.1}, and discuss
geometry behind that example. Associated to this example,
there is a family of $K3$ surfaces in the weighted projective $3$-space
$\PP^3[1,1,4,6]$ with weight $(q_1,q_2,q_3,q_4)=(1,1,4,6)$.  
There is a mirror pair of $K3$
surfaces $(S,\hat S)$. Here we know (cf. Belcastro \cite{B02})
that $$\mbox{Pic}(S)=U\quad\mbox{so that $\rho(S)=2$},$$
and that $S$ has a mirror
partner $\hat S$ whose Picard lattice is given by
$$\mbox{Pic}(\hat S)=U\perp (-E_8)^2\quad\mbox{ so that
$\rho(\hat S)=18$}.$$
The mirror $K3$ family can be defined by a hypersurface in the orbifold
ambient space $\PP^3[1,1,4,6]/G$ of degree $12$. Here $G$
is the discrete group of symmetry and can be given explicitly by 
$G=(\Z/3\Z)\times(\Z/2\Z)=\langle g_1\rangle\times\langle g_2\rangle$
where $g_1, g_2$ are 
generatoers whose actions are given by:
$$\begin{matrix} g_1 : (Y_1,Y_2,Y_3,Y_4) &\mapsto &(\zeta_3 Y_1, Y_2,
\zeta_3^{-1}Y_3,Y_4) \\
g_2 : (Y_1,Y_2,Y_3,Y_4) &\mapsto  &(Y_1,-Y_2, Y_3,-Y_4)\end{matrix}$$
(Here $\zeta_3=e^{2\pi i/3}$.)  The $G$-invariant monomials are
$$Y_1^{12},\,Y_2^{12},\,Y_3^3,\, Y_4^2,\, Y_1^6Y_2^6,\, Y_1Y_2Y_3Y_4.$$
The matrix of exponents is the following $6\times 5$ matrix
$$\begin{pmatrix} 12 & 0 & 0 & 0 & 1 \\
                 0 & 12 & 0 & 0 & 1 \\
                 0 & 0 & 3 & 0 & 1 \\
                 0 & 0 & 0 & 2 & 1 \\
                 6 & 6 & 0 & 0 & 1 \\
                 1 & 1 & 1 & 1 & 1 \end{pmatrix}$$
whose rank is $2$. Therefore we may conclude that the 
typical $G$-invariant polynomials is in $2$-parameters, and 
$\hat S$ can be defined by the following $2$-parameter family
of hypersurfaces of degree $12$ 
$$Y_1^{12}+Y_2^{12}+Y_3^3+Y_4^2+\lambda Y_1Y_2Y_3Y_4+\phi Y_1^6Y_2^6=0$$
in $\PP^3[1,1,4,6]/G$ with parameters $\lambda$ and $\phi$.

How do we conmpute the Picard--Fuchs differential equation of
this $K3$ family? 

Several physics articles are devoted to this question.
For instance, Klemm--Lerche--Mayr \cite{KLM}, 
Hosono--Klemm--Theisen--Yau \cite{HKTY}, Lian and Yau \cite{LY2}
determined the Picard--Fuchs differential equation of the
Calabi--Yau family using the GKZ hypergeometric system. Also it was 
noticed (cf. \cite{KLM}, \cite{LY2}) that the Picard--Fuchs system of
this family of $K3$ surfaces can be realized as the degeneration
of the Picard--Fuchs systems of the Calabi--Yau family.
The family of Calabi--Yau threefolds is a degree $24$ hypersurfaces
in $\PP^4[1,1,2,8,12]$ with $h^{1,1}=3$. The defining equation for this family is
given by
$$Z_1^{24}+Z_2^{24}+Z_3^{12}+Z_4^3+Z_5^2
-12\psi_0Z_1Z_2Z_3Z_4Z_5$$
$$-6\psi_1(Z_1Z_2Z_3)^6-\psi_2(Z_1Z_2)^{12}=0.$$
Its Picard--Fuchs system is given by
$$\begin{matrix} L_1&=&\Theta_x(\Theta_x-2\Theta_z)-12\,x(6\Theta_x+5)(6\Theta_x+1)\\
L_2&=& \Theta_y^2-y(2\Theta_y-\Theta_z+1)(2\Theta_y-\Theta_z)\\
L_3&=&\Theta_z(\Theta_z-2\Theta_y)-z(2\Theta_z-\Theta_x+1)(2\Theta_z-\Theta_x)
\end{matrix}
$$
where 
$$x=-\frac{2\psi_1}{1728^2\psi_0^6},\, y=\frac{1}{\psi_2^2}\quad\mbox{and}
\quad z=-\frac{\psi_2}{4\psi_1^2}$$ 
are deformation coordinates. 

Now the intersection of this Calabi--Yau hypersurface with the hyperplane
$Z_2-t\,Z_1=0$ gives rise to a family of $K3$ surfaces
$$b_0Y_1Y_2Y_3Y_4+b_1Y_1^{12}+b_2Y_2^{12}+b_3Y_3^3+b_4Y_4^2+b_5Y_1^6Y_2^6=0$$
in $\PP^3[1,1,4,6]$ of degree $12$. Taking $(b_0,b_1,b_2,b_3,b_4,b_5)=
(\lambda, 1,1,1,1,\phi)$ we obtain the $2$-parameter family of
$K3$ surfaces described above.
The Picard--Fuchs system of this $K3$ family is obtained by taking
the limit $y=0$ in the Picard--Fuchs system for the Calabi--Yau family:
$$\begin{matrix}
L_1&=&\Theta_x(\Theta_x-2\Theta_z)-12\,x(6\Theta_x+5)(6\Theta_x+1)\\
L_3&=&\Theta_z^2-z(2\Theta_z-\Theta_x+1)(2\Theta_z-\Theta_x)
\end{matrix}
$$
Further, if we intersect this $K3$ family with the hyperplane
$Y_2-s\,Y_1=0$, we obtain a family of elliptic curves:
$$c_0W_1W_2W_3+c_1W_1^6+c_2W_2^3+c_3W_3^2=0$$
in $\PP^2[1,2,3]$, whose Picard--Fuchs equation is given by
$$L=\Theta_x^2-12\,x(6\Theta_x+5)(6\Theta_x+1).$$ }
\end{ex}

Here we describe a relation of the Picard--Fuchs
system of the above family of $K3$ surfaces to the differential equation 
discussed in Example 4.1.

\begin{rem}\label{remark 6.1}
{\rm We note that, in view of our proof of Example 4.1, the process of
setting $z=0$ in the above Picard--Fuchs system $\{L_1,L_3\}$ is
equivalent to setting $t_1=0$ or $t_2=0$ in $x$ and $y$ in Example
4.1. Our Theorem 3.1 then implies that
$F(t)=(1-t)^{1/6}\,_2F_1(1/6,1/6;1;t)$ satisfies
$$
  (1+x)D_x^2F+xD_xF+\frac{5}{36}xF=0
$$
with $x=t/(1-t)$, or equivalently, (making a change of variable
$x\mapsto-x$)
$$
  x(1-x)F^{\prime\prime}+(1-2x)F^\prime-\frac5{36}F=0
$$
with $x=t/(t-1)$. That is,
$$
  (1-t)^{1/6}\,_2F_1\left(\frac16,\frac16;1;t\right)
 =\,_2F_1\left(\frac16,\frac56;1;\frac t{t-1}\right).
$$
This is the special case of the hypergeometric series identity
$$
  (1-t)^a\,_2F_1(a,b;c;t)=\,_2F_1\left(a,c-b;c;\frac t{t-1}\right).
$$
}
\end{rem}
\bigskip

We will discuss more examples of Picard--Fuchs systems of Calabi--Yau
threefolds and $K3$ surfaces,
which have already been considered by several people. For instance,
the articles \cite{HKTY}, \cite{HLY}, and \cite{KLM} obtained the Picard--Fuchs
operators for Calabi--Yau hypersurfaces with $h^{1,1}\leq 3$.  The next two 
examples consider Calabi--Yau hypersurfaces with $h^{1,1}>3$, and the paper of 
Lian and Yau \cite{LY2} addressed the question of determining the
Picard--Fuchs system of the families of $K3$ surfaces $\PP^3[1,1,2,2]$ of 
degree $6$ and $\PP^3[1,1,2,4]$ of degree $8$.  Their results are that 
\smallskip

(1) there is an elliptic fibration on these
$K3$ surfaces, and the Picard--Fuchs systems of the $K3$ families
can be derived from the Picard--Fuchs system of the elliptic
pencils, and that
\smallskip

(2) the solutions of the Picard--Fuchs systems for the
$K3$ families are given by ``squares'' of those for the
elliptic families.  
\smallskip

The system of partial differential equations considered
by Lian and Yau \cite{LY2} is
$$\begin{matrix}
L_1&=& \Theta_x(\Theta_x-2\Theta_z)-\lambda\,x(\Theta_x+\frac{1}{2}+\nu)(\Theta_x
+\frac{1}{2}-\nu) \\
L_2&=& \Theta_z^2-z(2\Theta_z-\Theta_x+1)(2\Theta_z-\Theta_x)
\end{matrix}
$$
and an ordinary differential equations
$$L=\Theta_x^2-\lambda\,x(\Theta_x+\frac{1}{2}+\nu)(\Theta_x+\frac{1}{2}-\nu)$$
where
$\Theta_x=x\frac{\partial}{\partial\,x}$, etc.) and $\lambda,\, \nu$
are complex numbers.

Also they noted that the $K3$ families correspond, respecitvely, to
the families 
of Calabi--Yau threefolds $\PP^4[1,1,2,4,4]$ of degree $12$ and 
$\PP^4[1,1,2,4,8]$ of degree $16$. However, the Picard--Fuchs
systems for the Calabi--Yau families are not explicitly
determined. 

\begin{ex}\label{example 6.2}
{\rm We now consider a family of $K3$ surfaces $\PP^3[1,1,2,4]$.
of degree $8$. This $K3$ family is realized as the degeneration
of the family of Calabi--Yau hypersurfaces $\PP^4[1,1,2,4,8]$ of
degree $16$ and $h^{1,1}=4$. The most generic defining
equation for this family is given by 
$$a_0Z_1Z_2Z_3Z_4Z_5+a_1Z_1^{16}+a_2Z_2^{16}+a_3Z_3^8+a_4Z_4^4+a_5Z_5^2
+a_6Z_3^2Z_4Z_5+a_7Z_1^8Z_2^8=0$$
Again the intersection with the hyperplane $Z_2-t\,Z_1=0$ gives rise
to a family of $K3$ surfaces $\PP^3[1,1,2,4]$:
$$Y_1^8+Y_2^8+Y_3^4+Y_4^2+\lambda Y_1Y_2Y_3Y_4+\phi Y_1^4Y_2^4=0$$
Let $S$ denote this family of $K3$ surfaces. Then 
$$Pic(S)=M_{(1,1),(1,1),0}\quad\mbox{with $\rho(S)=3$}.$$
The mirror family $\hat S$ exists and its Picard lattice
is 
$$Pic(\hat S)=E_8\perp D_7\perp U\quad\mbox{with $\rho(\hat S)=17$}.$$
The Picard lattices are determined by Belcastro \cite{B02}.
The intersection of this family of $K3$ surfaces with the
hyperplane $Y_2-s\,Y_2=0$ gives rise to the pencil of elliptic
curves
$$c_0W_1W_2W_3+c_1W_1^4+c_2W_2^4+c_3W_3^2=0$$
in $\PP^2[1,1,2]$ of degree $4$. This means that this family of
$K3$ surfaces has the elliptic fibration with section. 

Now translate this ``inductive'' structure to the Picard--Fuchs
systems.
The Picard--Fuchs system for the $K3$ family is given by
$$\begin{matrix}
L_1&=& \Theta_x(\Theta_x-2\Theta_z)-64\,x(\Theta_x+\frac{1}{2}+\frac{1}{4})
(\Theta_x+\frac{1}{2}-\frac{1}{4})\\
L_2&=& \Theta_z^2-z(2\Theta_z-\Theta_x+1)(2\Theta_z-\Theta_x)
\end{matrix}$$
and the Picard--Fuchs defferential equation of the elliptic
family is given by
$$L=\Theta_x^2-64\,x(\Theta_x+\frac{1}{2}+\frac{1}{4})(\Theta_x+\frac{1}{2}-\frac{1}{4})$$

The same remark as Remark 6.1 is valid for the Picard--Fuchs system $\{L_1,L_3\}$
which corresponds to Theorem 4.1 (b) with $a=1/3$.}
\end{ex}

\begin{ex}\label{example 6.3}
{\rm We consider a family of $K3$ surfaces $\PP^3[1,1,2,2]$ of
degree $6$. This $K3$ family is realized as the degeneration of
the family of Calabi--Yau hypersurfaces $\PP^4[1,1,2,4,4]$ of degree
$12$ and $h^{1,1}=5$:
$$a_0Z_1Z_2Z_3Z_4Z_5+a_1Z_1^{12}+a_2Z_2^{12}+a_3Z_3^6+a+4Z_4^3+a_5Z_5^3
+a_6Z_1^6Z_2^6=0.$$
The intersection of this Calabi--Yau hypersurface with the
hyperplane $Z_2-t\,Z_1=0$ gives rise to the family of $K3$ hypersurfaces
$\PP^3[1,1,2,4]$:
$$Y_1^6+Y_2^6+Y_3^3+Y_4^3+\lambda Y_1Y_2Y_3Y_4+\phi Y_1^3Y_2^3=0.$$
Let $S$ denote this family of $K3$ surfaces. Then
$$\mbox{Pic}(S)=M_{(1,1,1),(1,1,1),0}\quad\mbox{with $\rho(S)=4$}.$$
There is a mirror family of $K3$ surfaces, $\hat S$ with
$$\mbox{Pic}(\hat S)=E_8\perp D_4\perp A_2\perp U\quad\mbox
{with $\rho(\hat S)=16$}.$$
The Picard lattices are determined by Belcastro \cite{B02}.

The intersection of this $K3$ family with the hyperplane
$Y_2-s\,Y_1=0$ gives rise to the family of elliptic curves
$$c_0W_1W_2W_3+c_1W_1^3+c_2W_2^3+c_3W_3^3=0$$
in $\PP^2[1,1,1]$ of degree $3$.

The Picard--Fuchs system of this $K3$ family is
$$\begin{matrix}
L_1&=& \Theta_x(\Theta_x-2\Theta_z)-27\,x(\Theta_x+\frac{1}{2}+\frac{1}{6})(\Theta_x
+\frac{1}{2}-\frac{1}{6})\\
L_2&=& \Theta_z^2-z(2\Theta_z-\Theta_x+1)(2\Theta_z-\Theta_x)
\end{matrix}
$$
and the Picard--Fuchs differential equation for the elliptic family is
given by
$$L=\Theta_x^2-27\,x(\Theta_x+\frac{1}{2}+\frac{1}{6})(\Theta_x+\frac{1}{2}-\frac{1}{6})$$

We note that the same remark is valid for the Picard--Fuchs system $\{L_1,L_3\}$
corresponding to $a=1/4$ in Theorem 4.1(c).  }
\end{ex}

We will summarize the above discussions for the families of
$K3$ surfaces in the following form.

\begin{prop}\label{proposition 6.1}
{\sl The Picard--Fuchs systems
of families of $K3$ surfaces obtained by Lian and Yau \cite{LY2} can be
reconstructed starting from the modular forms (of two variables)
and then finding the
differential equations satisfied by them. In other words, the differential
equations satisfied by the modular forms (of two variables)
are realized as the Picard--Fuchs
differential equations of the families of $K3$ surfaces, establishing, 
in a sense, the ``modularity'' of the $K3$ families.} 
\end{prop}

\section{Picard--Fuchs differential equations of families of
$K3$ surfaces: Part II}
\label{sect7}

The purpose of this section is to study (one-parameter) families of $K3$ surfaces 
(some of which are realized as degenerations of some families of Calabi--Yau threefolds), 
whose mirror maps are expressed in terms of Hauptmodules for genus zero subgroups
of the form $\Gamma_0(N)^*$,  aiming to identify their Picard--Fuchs systems with 
differential equations assocaited to some to modular forms (of
two variable) (e.g., in Theorem 5.1). 

Dolgachev \cite{D96} has discussed several examples of families of $M_N$-polarized
$K3$ surfaces corresponding to $\Gamma_0(N)^*$ for small values of $N$, e.g., 
$N=1, 2$ and $3$.  

Lian and Yau \cite{LY1} have given examples of families of $K3$ surfaces and their Picard--Fuchs
differential equqtions of order $3$. The modular groups are genus zero subgroups of the form 
$\Gamma_0(N)^*$ where $N$ ranging from $1$ to $30$.  Here we try to analyze their examples and
their method in relation to our results in the section 5.

\begin{ex}\label{example 7.1}
{\rm We start with the hypergeometric equation:
$$t(1-t)f^{\prime\prime}+[1-(1+a+b)t]f^{\prime}-abf=0$$
in Theorem 5.1. Take $a=b=\frac{1}{4}$ and consider a one-parameter deformation of this
equation of the form:
$$t(1-t)f^{\prime\prime}+(1-\frac{3}{2}t)f^{\prime}-\frac{1}{16}(1-4\nu ^2)f=0$$
with a deformation parameter $\nu$.
This has a unique solution $f_0(t)$ near $t=0$ with $f_0(0)=1$, and a solution
$f_1(t)$ with $f_1(t)=f_0(t)\mbox{log}\,t+O(t)$. The inverse $t(q)$ of the
power series $q=\mbox{exp}(\frac{f_1(t)}{f_0(t)})=t + O(t^2)$ defines an invertible
holomorphic function in a disc, and $t(q)$ is the so-called mirror map. Put 
$$x(q)=\frac{1}{\lambda}t(\lambda q)\quad\mbox{for a given $\lambda$}.$$
One of the main results of Lian and Yau \cite{LY1} is
that for any complex numbers $\lambda,\,\nu$ with $\lambda\neq 0$, there is
a power series identity: 
$$_3F_2(\frac{1}{2},\frac{1}{2}+\nu,\frac{1}{2}-\nu;1,1;\lambda\,x(q))^2
=\frac{x^{\prime\,2}}{x^2(1-\lambda\,x)}$$
in the common domain of definitions of both sides. As before, $x^{\prime}(q)=D_q x(q)$.

For instance, take $(\lambda,\nu)=(2^63^3,\frac{1}{3}),\,(2^8,\frac{1}{4}),\,
(2^23^3,\frac{1}{6})$ and $(2^6,0)$, then these relations are given  
below. The mirror maps in these examples are expressed in terms of
Hauptmodules of genus zero modular groups of the form $\Gamma_0(N)^*$ ($\Gamma_0(1)^*=\Gamma$).

$$\begin{matrix}
\mbox{Label} & \mbox{Modular Relation} && & \mbox{Modular Group} \\ 
I &:\Bigg(\sum_{n=0}^{\infty}\frac{(6n)!}{(3n)!(n!)^3}\frac{1}{j(\tau)^n}\Bigg)^2& =& E_4(q) & \Gamma \\
II &:\Bigg(\sum_{n=0}^{\infty}\frac{(4n)!}{(n!)^4} x_2(\tau)^n\Bigg)^2& =
&\frac{x_2^{\prime\,2}}{x^2(1-256x)} & \Gamma_0(2)^*\\
III &:\Bigg(\sum_{n=0}^{\infty}\frac{(2n)!(3n)!}{(n!)^5} x_3(\tau)^n\Bigg)^2&=
&\frac{x_3^{\prime\,2}}{x_3^2(1-108x_3)} & \Gamma_0(3)^*\\
IV &:\Bigg(\sum_{n=0}^{\infty}\frac{(2n)!^3}{(n!)^6} x_4(\tau)^n\Bigg)^2&=
&\frac{x_4^{\prime\,2}}{x_4^2(1-64x_4)} & \Gamma_0(4)^*.
\end{matrix}
$$

\noindent Here $j(\tau),\, x_2(\tau), x_3(\tau)$ and $x_4(\tau)$ are
Hauptmodules for the genus zero subgroups $\Gamma,\,
\Gamma_0(2)^*,\,\Gamma_0(3)^*$ and $\Gamma_0(4)^*$, respectively.
Observe that in each modular relation, the right hand side is a modular
form of weight $4$ on the corresponding genus zero subgroup.

We know that $_3F_2(\frac{1}{2},\frac{1}{2}+\nu,\frac{1}{2}-\nu;1,1;\lambda\,x)$ is
a unique solution with the leading term $1+O(x)$ to the differential
operator
$$L=\Theta_x^3-\lambda\,x(\Theta_x+\frac{1}{2})(\Theta_x+\frac{1}{2}
+\nu)(\Theta_x+\frac{1}{2}-\nu).$$

In these examples, this differential operator is identified with
the Picard--Fuchs differential operator for a one-parameter family
of $K3$ surfaces, which are obtained by degenerating Calabi--Yau families.   
(Cf. Lian and Yau \cite{LY1}, Klemm, Lercher and Myer \cite{KLM}.)

$$\begin{matrix}
\quad & \mbox{CY family} & \mbox{$K3$ family} & \mbox{PF Operator}  \\
 I    & X(1,1,2,2,2)[8] & X(1,1,1,3)[6]      &  
\Theta^3-8x(6\Theta+5)(6\Theta+3)(6\Theta+1) \\
II    & X(1,1,2,2,6)[12] & X(1,1,1,1)[4]    & \Theta^3-4x(4\Theta+3)(4\Theta+2)(4\Theta+1)  \\ 
III   & X(1,1,2,2,2,2)[6,4] & X(1,1,1,1,1)[3,2] & \Theta^3-6x(2\Theta+1)(3\Theta+2)(3\Theta+1) \\
IV    & X(1,1,2,2,2,2,2)[4,4,4] & X(1,1,1,1,1,1)[2,2,2]& \Theta^3-8x(2\Theta+1)^3 
\end{matrix}
$$
\smallskip

The $K3$ families I and II have already been discussed in Lian--Yau 
\cite{LY2} (see also Verrill--Yui \cite{VY}) in relation to mirror maps.
The Picard group of I (resp. II) is given by
$$(-E_8)^2\oplus U_2\oplus <-4>\quad\text{(resp. $(-E_8)^2\oplus U_2\oplus <-2>$}).$$

The Calabi--Yau family III can be realized as a complete
intersection of the two hypersurfaces:
$$\begin{matrix}
Y_1^6+Y_2^6+Y_3^3+Y_4^3+Y_5^3+Y_6^3=0 \\
Y_1^4+Y_2^4+Y_3^2+Y_4^2+Y_5^2+Y_6^2=0
\end{matrix}
$$
This Calabi--Yau family has $h^{1,2}=68$ and $h^{1,1}=2$.
The $K3$ family is realized as the fiber space by setting
$$Y_1=Z_1^{1/2},\quad Y_2=\lambda Z_1^{1/2},\quad\text{and}\quad
Y_i=Z_i\quad\text{for $i=3,\cdots, 6$}$$
where $\lambda\in{\PP}^1$ is a parameter.
That is, we obtain a family of complete intersection $K3$ surfaces $X(1,1,1,1,1)[3,2]$:
$$\begin{matrix}
(1+\lambda^6)Z_1^3+Z_3^3+Z_4^3+Z_5^3+Z_6^3=0\\
(1+\lambda^4)Z_1^2+Z_3^2+Z_4^2+Z_5^2+Z_6^2=0
\end{matrix}
$$

{\bf Question: What is the Picard group of this K3 family?}
\smallskip

In the similar manner, the Calabi--Yau family IV can be realized as a
complete intersection of the three hypersurfaces:
$$\begin{matrix}
Y_1^4+Y_2^4+Y_3^2+Y_4^2+Y_5^2+Y_6^2+Y_7^2=0\\
Z_1^4+Z_2^4+Z_3^2+Z_4^2+Z_5^2+Z_6^2+Z_7^2=0\\
W_1^4+W_2^4+W_3^2+W_4^2+W_5^2+W_6^2+W_7^2=0
\end{matrix}$$
The $K3$ family is realized as the fiber space by setting
$$Y_1=Y_1^{\prime \frac{1}{2}},\, Y_2=\lambda Y_1^{\prime \frac{1}{2}}\quad\text{and}\quad
Y_i=Y_i^{\prime}\quad\text{for $i=3,\cdots, 7$}$$
and similarly for $Z_1,\, Z_2$ and $W_1,\, W_2$ where $\lambda\in\PP^1$ is a
parameter.

This gives rise to the $K3$ family $X(1,1,1,1,1,1)[2,2,2]$:
$$\begin{matrix}
(1+\lambda^4)Y_1^{\prime 2}+Y_3^{\prime 2}+Y_4^{\prime 2}+Y_5^{\prime 2}+Y_6^{\prime 2}
+Y_7^{\prime}=0\\
(1+\lambda^4)Z_1^{\prime 2}+Z_3^{\prime 2}+Z_4^{\prime 2}+Z_5^{\prime 2}+Z_6^{\prime 2}
+Z_7^{\prime 2}=0\\
(1+\lambda^4)W_1^{\prime 2}+W_3^{\prime 2}+W_4^{\prime 2}+W_5^{\prime 2}+W_6^{\prime 2}
+W_7^{\prime 2}=0
\end{matrix}
$$
\medskip

{\bf Question: What is the Picard group of this K3 family?}
\medskip
 
\noindent Here is the summary: 

(1) One starts with a Hauptmodule $x(=x(q))$ for a genus zero subgroup
$\Gamma_0(N)^*$; 

(2) then there associate a modular form $\frac{x^{\prime\,2}}{x\,r(x)}$ 
of weight $4$, 

(3) and a power series solution $\omega_0(x)$ of an order three differential operator;

(4) this differential operator coincides with the Picard--Fuchs differential
operator of a one-parameter family of $K3$ surfaces in weighted
projective spaces.}
\end{ex}

Lian and Yau \cite{LY1} further considered generalizations of the above
phenomenon, constructing many more examples.  Given a genus zero subgroup
of the form $\Gamma_0(N)^*$ and a Hauptmodul $x(q)$, constract
(by taking a Schwarzian derivative) a modular form $E$ of weight $4$ 
of the form $\frac{x^{\prime\,2}}{x\,r(x)}$
and a differential operator $L$ whose monodromy has maximal unipotency
at $x=0$, such that $L\,E^{1/2}=0$.  Further, identify $L$ as the Picard--Fuchs
differential operator of a family of $K3$ surfaces. Let $\omega_0(x)$ denotes
the fundamental period of this manifold. Then it should be subject to
the modular relation
$$\omega_0(x)^2=\frac{x^{\prime\,2}}{x\,r(x)}$$ 
\medskip

How do we associate modular forms of weight $(1,1)$ corresponding to the 
groups $\Gamma_0(N)^*\times \Gamma_0(N)^*$ in this situation?

Taking the square root of both sides of the modular relation, we  
obtain that $\omega_0(x)^{1/2}$ is a modular form (of one
variable) of weight $1$ for the group $\Gamma_0(N)^*$. 
Then taking $\omega_0(q_1)\omega_0(q_2)$, we see that this is a modular form
for $\Gamma_0(N)^*\times \Gamma_0(N)^*$ of weight $(1,1)$. 
Then this modular form (of two variables) satisfies a differential equation, 
which may be identified with the Picard--Fuchs differential equation of
the K3 family considered above. We summarize the above discussion
in the following proposition.

\begin{prop}\label{proposition 7.1}
{\sl The examples I--IV above are related to our Theorem 5.2. Indeed, the 
connection is established by the identity
  $$
    _2F_1\left(a,b;a+b+\frac12;z\right)^2=\,
    _3F_2\left(2a,a+b,2b;a+b+\frac12,2a+2b;z\right).
  $$
  More explicitly, the examples I--IV correspond to the cases $(1/12,5/12)$,
  $(1/8,3/8)$, $(1/6,1/3)$, and $(1/4,1/4)$, respectively.}
\end{prop}

Note that the generalized hypergeometric series
  $_3F_2(\alpha_1,\alpha_2,\alpha_3;1,1;z)$ satisfies the differential
  equation of the form: 
$$[\Theta_z^3-\lambda\,z(\Theta_z+\alpha_1)(\Theta_z+\alpha_2)(\Theta_z+\alpha_3)]f=0$$ 
for some $\alpha_1,\alpha_2,\alpha_3\in{\Q}$ and $\lambda\in{\Q},\, \neq 0$.

A natural question we may ask now is: Is is possible to construct families of
$K3$ surfaces corresponding to Theorem 5.2 from this observation?   

When the order $3$ differential equation of this form becomes the symmetric
square of an order $2$ differential equation, and if the order $2$
differential equation is realized as the Picard--Fuchs differential
equation of a family of elliptic curves, we may be able to construct
a family of $K3$ surfaces using the method of Long \cite{Long}, especially
when the Picard number of the $K3$ family in question is $19$ or $20$.
In fact, Rodriguez--Villegas \cite{RV03} has discussed $4$ families
of $K3$ surfaces which fall into this class.

However, at the moment, we do not know if there are readily available methods 
for constructing $K3$ families starting from differential equations.

\begin{rem}
{\rm If we consider the order $4$ generalized hypergeomtric series, there
are $14$ differential equations are of the form
$$[\Theta_z^4-\lambda\,z(\Theta_z+\alpha_1)(\Theta_z+\alpha_2)(\Theta_z+\alpha_3)(\Theta_z+\alpha_4)]f
=0$$
for some $\alpha_i\in{\Q}$ and $\lambda\in {\Q},\, \neq 0$.
These $14$ differential operators can be found in Almkvist--Zudilin \cite{AZ}.
As Doran and Morgan \cite{DorMor05} explained, only 13 of the 14 such operators
are known to be realizable as the Picard--Fuch differential operator for
a family of smooth Calabi--Yau threefolds with $h^{2,1}=1$. 
For the $13$ cases, Klemm and Theisen \cite{KT} (see also Villegas \cite{RV03})
found the corresponding families of Calabi--Yau 
threefolds in weighted projective spaces.  
The missing case is $(\alpha_1,\alpha_2,\alpha_3,\alpha_4)=(1/12, 5/12, 7/12, 11/12)$.
For more thorough discussions on this topic, the reader should consult
the article of Doran and Morgan \cite{DorMor05}, where classify integral
monodromy representations. } 
\end{rem}

\section{Generalizations and open problems}
\label{sect8}

\begin{pr}\label{problem 8.1}
{\rm We have determined differential equations satisfied by modular forms
(of two variables) of weight $(1,1)$. The arguments can be generalized to 
modular forms (of two variables) of any weight $(k_1,k_2)$, using the 
result of Yang \cite{Yg04}.
However, differerntial equations satisfied by them are getting 
too big to display.}
\end{pr}

\begin{pr}\label{problem 8.2}
{\rm A natural generalization is to consider modular forms of
three (or more than three) variables 
$F(\tau_1,\tau_2,\tau_3)$ of weight $(k_1,k_2,k_3)$ on $\Gamma_1\times
\Gamma_2\times \Gamma_3$.   

Examples of this kind should correspond to Picard--Fuchs
differential equations of families of Calabi--Yau threefolds, or Picard--Fuchs
differential equations of degenerate families of Calabi--Yau fourfolds.}
\end{pr}

\section*{Acknowledgments}

The collaboration for this work started when Y. Yang visited N. Yui
at Queen's University, Kingston Canada, in August 2004. The work was 
completed at Tsuda College, Tokyo Japan 2005,
and revisions were incorporated at Tsuda College again in the summer
of 2006. Both authors were visiting professors at that institution
in June 2005, and August 2006. We thank the hospitality of Tsuda College. 

We thank Jan Stienstra for his interest and helpful comments and
suggestions, and Chuck Doran and Don Zagier for their comments
on the earlier version(s) of this paper. 

\vskip 1.5cm

%%%%%%%%%%%%%%%%%%%%%%%%%%%%%%%%%%%%%%%%%%%%%%%%%%%%%%%%%%%%%%%%%%%%%%%%%%%

\end{document}